%

\catcode`\@=11

%
%
\def\bibn@me{R\'ef\'erences}
\def\bibliographym@rk{\centerline{{\sc\bibn@me}}
	\sectionmark\section{\ignorespaces}{\unskip\bibn@me}
	\bigbreak\bgroup
	\ifx\ninepoint\undefined\relax\else\ninepoint\fi}
%
%
%
\let\refsp@ce=\
\let\bibleftm@rk=[
\let\bibrightm@rk=]
%
%
%
\def\numero{n\raise.82ex\hbox{$\fam0\scriptscriptstyle
o$}~\ignorespaces}
%
%
\newcount\equationc@unt
\newcount\bibc@unt
\newif\ifref@changes\ref@changesfalse
\newif\ifpageref@changes\ref@changesfalse
\newif\ifbib@changes\bib@changesfalse
\newif\ifref@undefined\ref@undefinedfalse
\newif\ifpageref@undefined\ref@undefinedfalse
\newif\ifbib@undefined\bib@undefinedfalse
\newwrite\@auxout
%
%
%
%
%
%
%
%
\def\re@dreferences#1#2{{%
	\re@dreferenceslist{#1}#2,\undefined\@@}}
\def\re@dreferenceslist#1#2,#3\@@{\def\next{#2}%
	\expandafter\ifx\csname#1@@\meaning\next\endcsname\relax
	??\immediate\write16
	{Warning, #1-reference "\next" on page \the\pageno\space
	is undefined.}%
	\global\csname#1@undefinedtrue\endcsname
	\else\csname#1@@\meaning\next\endcsname\fi
	\ifx#3\undefined\relax
	\else,\refsp@ce\re@dreferenceslist{#1}#3\@@\fi}
%
%
%
\def\newlabel#1#2{{\def\next{#1}\newl@bel#2}}
\def\newl@bel#1#2{%
	\expandafter\xdef\csname ref@@\meaning\next\endcsname{#1}%
	\expandafter\xdef\csname pageref@@\meaning\next\endcsname{#2}}
\def\label#1{{%
	\toks0={#1}\message{ref(\lastref) \the\toks0,}%
	\ignorespaces\immediate\write\@auxout%
	{\noexpand\newlabel{\the\toks0}{{\lastref}{\the\pageno}}}%
	\def\next{#1}%
	\expandafter\ifx\csname ref@@\meaning\next\endcsname\lastref%
	\else\global\ref@changestrue\fi%
	\newlabel{#1}{{\lastref}{\the\pageno}}}}
\def\ref#1{\re@dreferences{ref}{#1}}
\def\pageref#1{\re@dreferences{pageref}{#1}}
%
%
\def\bibcite#1#2{{\def\next{#1}%
	\expandafter\xdef\csname bib@@\meaning\next\endcsname{#2}}}
\def\cite#1{\bibleftm@rk\re@dreferences{bib}{#1}\bibrightm@rk}
%
%
\def\beginthebibliography#1{\bibliographym@rk
	\setbox0\hbox{\bibleftm@rk#1\bibrightm@rk\enspace}
	\parindent=\wd0
	\global\bibc@unt=0
	\def\bibitem##1{\global\advance\bibc@unt by 1
		\edef\lastref{\number\bibc@unt}
		{\toks0={##1}
		\message{bib[\lastref] \the\toks0,}%
		\immediate\write\@auxout
		{\noexpand\bibcite{\the\toks0}{\lastref}}}
		\def\next{##1}%
		\expandafter\ifx
		\csname bib@@\meaning\next\endcsname\lastref
		\else\global\bib@changestrue\fi%
		\bibcite{##1}{\lastref}
		\medbreak
		\item{\hfill\bibleftm@rk\lastref\bibrightm@rk}%
		}
	}
\def\endthebibliography{\egroup\par}
%
%
    \outer\def\bye{\@closeaux
    	\par\vfill\supereject\end}
%
\def\@closeaux{\closeout\@auxout
	\ifref@changes\immediate\write16
	{Warning, changes in references.}\fi
	\ifpageref@changes\immediate\write16
	{Warning, changes in page references.}\fi
	\ifbib@changes\immediate\write16
	{Warning, changes in bibliography.}\fi
	\ifref@undefined\immediate\write16
	{Warning, references undefined.}\fi
	\ifpageref@undefined\immediate\write16
	{Warning, page references undefined.}\fi
	\ifbib@undefined\immediate\write16
	{Warning, citations undefined.}\fi}
%
%
\immediate\openin\@auxout=\jobname.aux
\ifeof\@auxout \immediate\write16
     {Creating file \jobname.aux}
\immediate\closein\@auxout
\immediate\openout\@auxout=\jobname.aux
\immediate\write\@auxout {\relax}%
\immediate\closeout\@auxout
\else\immediate\closein\@auxout\fi
%
%
\input\jobname.aux \par
\immediate\openout\@auxout=\jobname.aux
%
%

\def\bibn@me{R\'ef\'erences bibliographiques}
\catcode`@=11
\def\bibliographym@rk{\bgroup}
%
%
\outer\def\bye{ 	\par\vfill\supereject\end}

\magnification=1200

\font\tenbfit=cmbxti10 
\font\sevenbfit=cmbxti10 at 7pt 
\font\sixbfit=cmbxti5 at 6pt 

\newfam\mathboldit 

\textfont\mathboldit=\tenbfit
  \scriptfont\mathboldit=\sevenbfit
   \scriptscriptfont\mathboldit=\sixbfit

\def\bfit           
{\tenbfit           
   \fam\mathboldit 
}

\def\QQ{{\bf {Q}}}

\def\ZZ{{\bf {Z}}}
\def\CC{{\bf {C}}}
\def\RR{{\bf {R}}}

\def\Q{{\bf {Q}}}
\def\K{{\bf {K}}}

\def\R{{\bf R}}

\def\Bad{{\bfit Bad}}

\def\house#1{\setbox1=\hbox{$\,#1\,$}%
\dimen1=\ht1 \advance\dimen1 by 2pt \dimen2=\dp1 \advance\dimen2 by 2pt
\setbox1=\hbox{\vrule height\dimen1 depth\dimen2\box1\vrule}%
\setbox1=\vbox{\hrule\box1}%
\advance\dimen1 by .4pt \ht1=\dimen1
\advance\dimen2 by .4pt \dp1=\dimen2 \box1\relax}

\def\Norm{{\rm Norm}}

  \def\eps{{\varepsilon}}

\def\sm{\smallskip}  

\def\build#1_#2^#3{\mathrel{\mathop{\kern 0pt#1}\limits_{#2}^{#3}}}

\def\date {le\ {\the\day}\ \ifcase\month\or janvier
\or fevrier\or mars\or avril\or mai\or juin\or juillet\or
ao\^ut\or septembre\or octobre\or novembre
\or d\'ecembre\fi\ {\oldstyle\the\year}}

\font\fivegoth=eufm5 \font\sevengoth=eufm7 \font\tengoth=eufm10

\newfam\gothfam \scriptscriptfont\gothfam=\fivegoth
\textfont\gothfam=\tengoth \scriptfont\gothfam=\sevengoth

\def\pro{\noindent {\it Proof. }}

\def\smallsquare{\vbox{\hrule\hbox{\vrule height 1 ex\kern 1 ex\vrule}\hrule}}
\def\cqfd{\hfill \smallsquare\vskip 3mm}

\def\Q{{\bf {Q}}}  
\def\K{{\bf {K}}} 



\centerline{}

\vskip 4mm

\centerline{
\bf On a mixed problem in Diophantine 
approximation}

\vskip 8mm
\centerline{Yann B{\sevenrm UGEAUD}
\ \& \ Bernard de M{\sevenrm ATHAN}
\footnote{}{\rm 
}
}

{\narrower\narrower
\vskip 12mm

\proclaim Abstract. {
Let $d$ be a positive integer.
Let $p$ be a prime number.
Let $\alpha$ be a real algebraic number
of degree $d+1$. We establish that
there exist a positive constant $c$ 
and infinitely many algebraic numbers $\xi$ 
of degree $d$ such that
$|\alpha - \xi| \cdot \min\{|\Norm(\xi)|_p,1\} < 
c H(\xi)^{-d-1} \, (\log 3 H(\xi))^{-1/d}$. 
Here, $H(\xi)$ and $\Norm(\xi)$ denote the
na\"\i ve height of $\xi$ and
its norm, respectively. This extends 
an earlier result of de Mathan and Teuli\'e
that deals with the case $d=1$.
}

}

\vskip 6mm
\vskip 15mm

\centerline{\bf 1. Introduction}

\vskip 6mm

In analogy with the Littlewood conjecture,
de Mathan and Teuli\'e \cite{BdMTe} 
proposed recently a `mixed Littlewood
conjecture'.
For any prime number $p$, the usual
$p$-adic absolute value $| \cdot |_p$ is
normalized in such a way that $|p|_p = p^{-1}$.
We denote by $\Vert \cdot \Vert$ the distance
to the nearest integer.

\proclaim De Mathan--Teuli\'e conjecture.
For every real number $\alpha$ and every prime
number $p$, we have
$$
\inf_{q \ge 1} \, q \cdot \Vert q \alpha \Vert \cdot
\vert q \vert_p = 0. \eqno (1.1)
$$

Obviously, the above conjecture holds if $\alpha$ is rational or
has unbounded partial quotients in its continued fraction expansion. 
Thus, it only remains
to consider the case
when $\alpha$ is an element of the set $\Bad_1$ of badly approximable
real numbers, where
$$
\Bad_1  = \{ \alpha \in \R : \inf_{q \ge 1} \,
q \cdot \Vert q \alpha \Vert > 0\}.
$$
De Mathan and Teuli\'e \cite{BdMTe} proved that (1.1) holds
for every quadratic real number $\alpha$ (recall
that such a number is in $\Bad_1$) but,
despite several recent results \cite{EiKl07,BDM}, 
the general conjecture is still unsolved.

If we rewrite (1.1) under the form 
$$
\inf_{a, q\ge1, \, {\rm gcd}(a, q) = 1} \, 
q^2 \cdot \biggl|\alpha-{a\over q}\biggr| \cdot |q|_p=0\,,
$$ 
then we have $|q|_p=\min\{|\Norm(q/a)|_p,1\}$. 
Hence, upon replacing $\alpha$ by $1/\alpha$, 
the de Mathan--Teuli\'e conjecture can be reformulated as follows: 
For every irrational real number $\alpha$, for every prime
number $p$ and every positive real
number $\eps$, there exists a non-zero rational number $\xi$ satisfying
$$
|\alpha - \xi| \cdot \min\{|\Norm(\xi)|_p,1\} < \eps H(\xi)^{-2}.
$$
Throughout this paper, the height $H(P)$ of an integer
polynomial $P(X)$ is the maximal of the absolute values
of its coefficients. The height $H(\xi)$ of an
algebraic number $\xi$ is the height 
of its minimal defining polynomial over the rational integers
$a_0 + a_1 X + \ldots + a_d X^d$, and the norm
of $\xi$, denoted by $\Norm(\xi)$, is the rational
number $(-1)^d a_0 / a_d$.

This reformulation suggests us to ask the following question.

\proclaim Problem 1. 
Let $d$ be a positive integer. Let
$\alpha$ be a real number that is not algebraic of
degree less than or equal to $d$. 
For every prime
number $p$ and every positive real
number $\eps$, does there exist a non-zero real
algebraic number $\xi$ of degree at most $d$ satisfying
$$
|\alpha - \xi| \cdot \min\{|\Norm(\xi)|_p,1\} < \eps H(\xi)^{-d-1} ?
$$

The answer to Problem 1 is clearly positive,
unless (perhaps) when
$\alpha$ is an element of the set $\Bad_d$ of 
real numbers that are badly approximable by algebraic numbers of degree 
at most $d$, where
$$
\eqalign{
\Bad_d  = & \{ \alpha \in \R : \hbox{There exists $c > 0$ such that 
$|\alpha - \xi| > c H(\xi)^{-d-1}$}, \cr
& \hbox{for all algebraic numbers $\xi$ of degree
at most $d$} \}. \cr}
$$
For $d \ge 1$, the set $\Bad_d$ contains the set
of algebraic numbers of degree $d+1$, but it remains 
an open problem to decide whether this 
inclusion is strict for $d \ge 2$; 
see the monograph~\cite{BuLiv} for more information.
The purpose of the present note is to 
give a positive answer to Problem 1
for every positive integer $d$ and
every real algebraic number $\alpha$ of degree $d+1$. This extends
the result from \cite{BdMTe} which deals with the case $d=1$.

\vskip 8mm

\centerline{\bf 2. Results}

\vskip 6mm

Throughout this paper, for a prime 
number $p$, a number field $\K$,
and a non-Archime\-dean place $v$ on $\K$ lying above $p$, we normalize
the absolute value $| \cdot |_v$ in such a way that
$| \cdot |_v$ and $|\cdot |_p$ coincide on $\Q$.

Our main result includes
a positive answer to Problem 1 when 
$\alpha$ is a real algebraic number of degree $d+1$.

\proclaim Theorem 1.
Let $d$ be a positive integer.
Let $\alpha$ be a real algebraic number
of degree $d+1$ and denote by $r$ the unit rank of $\QQ(\alpha)$.
Let $p$ be a prime number.
There exist positive constants $c_1, c_2, c_3$ 
and infinitely many real algebraic numbers $\xi$ 
of degree $d$ such that
$$
|\alpha -\xi| < c_1 H(\xi)^{-d-1}, \eqno (2.1)
$$
$$
|\xi|_v < c_2 (\log 3 H(\xi))^{-1/(rd)}, \eqno (2.2)
$$
for every absolute value $| \cdot |_v$ on $\QQ(\xi)$
above the prime $p$, and
$$
|\alpha - \xi| \cdot \min\{|\Norm(\xi)|_p,1\} < c_3 H(\xi)^{-d-1}
\, (\log 3 H(\xi))^{-1/r}.  \eqno (2.3)
$$

Theorem 1 extends Th\'eor\`eme 2.1
of \cite{BdMTe} that is only concerned
with the case $d=1$.

Under the assumptions of Theorem 1,
Wirsing \cite{Wir} established
that there are infinitely many real algebraic numbers $\xi$
satisfying (2.1). 

The proof of Theorem 1 is very much inspired by a paper
of Peck \cite{Peck} on simultaneous rational approximation
to real algebraic numbers. Roughly speaking, we use a method
dual to Peck's to construct integer polynomials $P(X)$
that take small values at $\alpha$, and we need
an extra argument to ensure that our polynomials
have a root $\xi$ very close to $\alpha$.

De Mathan \cite{BdM05} used the theory of linear forms
in non-Archimedean logarithms to prove that Theorem 1 for $d=1$
is best possible, in the sense that the absolute value 
of the exponent 
of $(\log 3 H(\xi))$ in (2.2) cannot be too large.
Next theorem extends this result to all values of $d$.

\proclaim Theorem 2.
Let $p$ be a prime number, $d$ a positive integer
and $\alpha$ a real algebraic number
of degree $d+1$. Let $\lambda$ be a positive
real number. There exists a
positive real number $\kappa=\kappa(\lambda)$ 
such that 
for every non-zero real algebraic number $\xi$ 
of degree $d$ satisfying 
$$
|\alpha -\xi|\le\lambda H(\xi)^{-d-1} \eqno (2.4)
$$ 
we have 
$$
|\xi|_v\ge(\log 3H(\xi))^{-\kappa}\,
$$ 
for at least one absolute value $| \cdot |_v$ on 
$\QQ(\xi)$ above the prime $p$.


As in \cite{BdM05}, the proof of
Theorem 2 rests on the theory of linear forms
in non-Archime\-dean logarithms.

Let $d$ be a positive integer.
We recall that it follows from the $p$-adic
version of the Schmidt Subspace Theorem that
for every algebraic number $\alpha$
of degree $d+1$ and for every positive real number
$\eps$, there are only finitely many non-zero
integer polynomials $P(X)=a_0+a_1X+ \ldots +a_dX^d$ 
of degree at most $d$, with $a_0\not=0$, 
that satisfy 
$$
|P(\alpha)| \cdot |a_0|_p < H(P)^{-d-\eps}\,.
$$ 
Let $\xi$ be a real algebraic number of degree at most $d$, 
and denote by $P(X)= a_0+a_1X+ \ldots +a_dX^d$ 
its minimal defining polynomial
over $\ZZ$. Then,
$$
\min\{|{\rm Norm} (\xi)|_p,1\}\ge|a_0|_p\,
$$  
and there exists a constant $c(\alpha)$, depending
only on $\alpha$, such that
$$
|P(\alpha)|\le c(\alpha) \,  H(\xi) \cdot |\xi-\alpha|.
$$ 
Let $\eps$ be a positive real number.
Applying the above statement deduced from the $p$-adic
version of the Schmidt Subspace Theorem to these
polynomials $P(X)$, we deduce that 
$$
|\alpha -\xi| \cdot \min\{|{\rm Norm}(\xi)|_p,1\}
\ge H(P)^{-d-1-\eps}\,
$$ 
holds if $H(P)$ is sufficiently large.
This implies that if $\xi$ satisfies (2.4)
and if $H(\xi)$ is sufficiently large, then
we have 
$$
|{\rm Norm}(\xi)|_p \ge H(\xi)^{-\eps}\,,
$$ 
accordingly 
$$
\max_{v|p}|\xi|_v \ge H(\xi)^{-\eps/d}\,.
$$ 
The result of Theorem 2 is 
more precise, however we cannot 
obtain a good lower bound for $|{\rm Norm}(\xi)|_p$.

We conclude this section by pointing out
that Einsiedler and Kleinbock \cite{EiKl07} 
showed that a slight modification
of the de Mathan--Teuli\'e conjecture easily follows from
a theorem of Furstenberg \cite{Fur67,Bos94}.

\proclaim Theorem EK.
Let $p_1$ and $p_2$ be distinct prime numbers.
Then 
$$
\inf_{q \ge 1} \, q
\cdot  \Vert q\alpha  \Vert \cdot \vert q\vert_{p_1} 
\cdot \vert q\vert_{p_2} =0 
$$
holds for every real number $\alpha$.

In view of Theorem EK, we formulate the 
following question, presumably easier to solve
than Problem 1.

\proclaim Problem 2. 
Let $d$ be a positive integer. Let
$\alpha$ be a real number that is not algebraic of
degree less than or equal to $d$. 
For every distinct prime
numbers $p_1$, $p_2$ and every positive real
number $\eps$, does there exist a non-zero real
algebraic number $\xi$ of degree at most $d$ satisfying
$$
|\alpha - \xi| \cdot \min\{|\Norm(\xi)|_{p_1}, 1\} \cdot
\min\{|\Norm(\xi)|_{p_2}, 1\} 
< \eps H(\xi)^{-d-1} ?
$$

Theorem EK gives a positive answer to Problem 2 when $d=1$.

The sequel of the paper is organized as follows.
We gather several auxiliary results in Section 3,
and Theorems 1 and 2 are established in Sections
4 and 5, respectively.

In the next sections, we fix a real algebraic number field $\K$ of
degree $d+1$. The notation $A \ll B$ means, unless
specific indications, that the
implicit constant depends on $\K$. Furthermore,
we write $A \asymp B$ if we have simultaneously 
$A \ll B$ and $B \ll A$.

\vskip 6mm
\goodbreak

\centerline{\bf 3. Auxiliary lemmas}

\vskip 4mm

Let $\K$ be a real algebraic number field of degree $d+1$.
Let ${\cal O}$ denote its ring of integers, and let
$\alpha_0=1, \alpha_1, \ldots ,\alpha_d$ be a basis of $\K$. 
Let $D$ be a positive integer satisfying
$$
D(\ZZ+\alpha_1\ZZ+ \ldots +\alpha_d\ZZ)\subset{\cal O}\subset 
{1\over D}(\ZZ+\alpha_1\ZZ+ \ldots +\alpha_d\ZZ)
$$ 
and the corresponding inequalities for the dual
basis $\beta_0, \ldots ,\beta_d$ defined by
$$
{\rm Tr}(\alpha_i\beta_j)=\delta_{i,j},
$$
where ${\rm Tr}$ is the trace and $\delta_{i, j}$ is the
Kronecker symbol.

We denote by $\sigma_0={\rm Id},  \ldots ., \sigma_d$, the complex
embeddings of $\K$, numbered in such a way that 
$\sigma_0, \ldots , \sigma_{r_1 - 1}$ are real,
$\sigma_{r_1}, \ldots , \sigma_d$ are imaginary
and $\sigma_{r_1+r_2+j}=\overline{\sigma_{r_1+j}}$ for $0\le j<r_2$.
Set also $r=r_1+r_2 - 1$, and let $\eps_1, \ldots , \eps_r$ 
be multiplicatively independent
units in $\K$.

\proclaim Lemma 1.
Let $\eta$ be a unit in ${\cal O}$ such that
$-1 < \eta < 1$ and define
the real number $N$ by $|\eta|=N^{-1}$.
The conditions
$$
|\sigma_j(\eta)|\asymp N^{1/d}\,,{\hskip2cm}0<j\le d\,,\eqno (3.1)
$$
and
$$
|\sigma_i(\eta)|\asymp|\sigma_j(\eta)|\,,{\hskip2cm}
0<i<j\le d\,,\eqno(3.2)
$$
are equivalent. Let $\gamma\not=0$ be in $\K$ and let $\Delta$
be a positive integer such that $\Delta\gamma\in {\cal O}$. 
If $\eta$ satisfies (3.1) or (3.2), write
$$
\gamma\eta=a_0+ \ldots +a_d\alpha_d\,,
$$
with $a_0, \ldots , a_d$ in $\QQ$.
We have $D\Delta a_k\in\ZZ$ for $k=0,\ldots,d$ and
$$
\max_{k=0,\ldots,d}|a_k|\asymp N^{1/d},
$$ 
where the implicit constants depend on $\gamma$.

\pro
Since $\eta$ is a unit, we have
$$
\prod_{0\le j\le d}\sigma_j(\eta)=\pm1\,.
$$
and (3.1) and (3.2) are clearly equivalent.
The formula
$$
a_k={\rm Tr}(\gamma\eta\beta_k)=\gamma\eta\beta_k+
\sum_{j=1}^d\sigma_j(\eta)\sigma_j(\gamma\beta_k)
$$
implies that if $\eta$ satisfies (3.1), then
$$
|a_k|\ll N^{1/d}\,,{\hskip2cm}0\le k\le d\,.
$$
Combined with
$$
\sigma_1(\gamma)\sigma_1(\eta)=a_0+\ldots+a_d\sigma_1(\alpha_d)\,,
$$ 
this shows that
$$
N^{1/d}\asymp|\sigma_1(\eta)|\ll\max_{k=0,\ldots,d}|a_k|.
$$
The proof of the lemma is complete.
\cqfd

Let $\alpha$ be a real algebraic number of degree $d+1$. 
We keep the above notation with the field $\K=\QQ(\alpha)$
and the basis $1,\alpha,\ldots,\alpha^d$ of $\K$ over $\QQ$,
and we display an immediate consequence of Lemma 1.

\proclaim Corollary 1.
Let $\eta$ be a unit in ${\cal O}$ such that
$-1 < \eta < 1$ and set $N = |\eta|^{-1}$.
Then 
$$
D\Delta\gamma\eta=P(\alpha),
$$ 
where $P(X)$ is an integral polynomial of 
degree at most $d$ satisfying
$$
H(P)\asymp N^{1/d}, \quad |P(\alpha)|\asymp N^{-1}\,,
$$ 
and thus
$$
|P(\alpha)|\asymp H(P)^{-d}.
$$

Denote by $\tau_j$, $j = 0, \ldots , d$ the embeddings
of $\K$ into $\CC_p$. Recall that the absolute
value $| \cdot |_p$ on $\QQ$ has a unique extension to
$\CC_p$, that we also denote by $| \cdot |_p$.
In Lemmata 2 to 4 below we work in $\CC_p$.
Let $P(X)$ be an irreducible integer polynomial
of degree $n \ge 1$. Let $\xi$ be a complex root of $P(X)$
and $\xi_1, \ldots , \xi_n$ be the roots of $P(X)$
in $\CC_p$. We point out that the sets
$$
\{ |\xi|_v : \hbox{$v$ is above $p$ on $\QQ(\xi)$} \}
$$
and
$$
\{ |\xi_i|_p : 1 \le i \le n \}
$$
coincide, since all the absolute values $| \cdot |_v$
and $| \cdot |_p$ coincide on $\QQ$.

Keeping the notation of Lemma 1, 
we have the following auxiliary result.

\proclaim Lemma 2.
Assume that $\gamma=\alpha_d$. Then
$$
|a_k|_p\ll\max_{0\le j\le d}|\tau_j(\eta)-1|_p\,, {\hskip2cm} 0\le k<d\,,
$$ 
and
$$ 
|a_d-1|_p\ll\max_{0\le j\le d}|\tau_j(\eta)-1|_p\,.
$$

\pro
Since
$$
{\rm Tr}(\alpha_d\beta_k)=0\,,
\quad \hbox{for $k=0, \ldots , d-1$},
$$ 
we get
$$
a_k={\rm Tr}(\gamma\eta\beta_k)={\rm Tr}(\alpha_d(\eta-1)\beta_k)=
\sum_{j=0}^d(\tau_j(\eta)-1) \tau_j(\alpha_d\beta_k),
$$ 
and deduce that
$$
|a_k|_p\ll\max_{0\le j\le d}|\tau_j(\eta)-1|_p\,, {\hskip2cm} 0\le k<d.
$$
It follows from
$$
{\rm Tr}(\alpha_d\beta_d)=1\,
$$ 
that
$$
a_d=1+{\rm Tr}(\alpha_d\beta_d(\eta-1))
=1+\sum_{j=0}^d(\tau_j(\eta)-1) \tau_j(\alpha_d\beta_d),
$$ 
and we derive that
$$ 
|a_d-1|_p\ll\max_{0\le j\le d}|\tau_j(\eta)-1|_p\,.
$$
This concludes the proof.
\cqfd

\proclaim Lemma 3.
Let $0<\delta<1$. 
There exist arbitrarily large positive real numbers $H$ 
and units $\eta$ satisfying $\eta=H^{-d}$,
$$
\left\vert{\sigma_j(\eta)\over\sigma_1(\eta)}-
1\right\vert\le\delta, {\hskip1cm} 2\le j\le d, \eqno (3.3)
$$
and
$$
|\tau_j (\eta)-1|_p\ll(\log H)^{-1/r}, {\hskip1cm} 0\le j \le d.
$$ 

\pro
By taking suitable powers of the 
units $\eps_1, \ldots , \eps_r$, we can assume that they
are all positive, as well as their real conjugates, 
and that $|\tau_j(\eps_i)-1|_p<p^{-1/(p-1)}$ 
for $i=1, \ldots , r$ and $j = 0, \ldots , d$.
This is possible since $|\tau_j(\eps_i)|_p=1$ for $i=1, \ldots , r$
and $j = 0, \ldots , d$. This allows us to consider the $p$-adic 
logarithms of each $\tau_j(\eps_i)$.
Our aim is to construct a suitable unit $\eta$ of the form
$$
\eta=\eps_1^{\mu_1p^s}\ldots\eps_r^{\mu_rp^s},
$$ 
where $\mu_i\in\ZZ$. The conditions for (3.3) are then
$$
p^s\left\vert\mu_1\log{|\sigma_j(\eps_1)|\over|\sigma_1(\eps_1)|}+
\ldots+\mu_r\log{|\sigma_j(\eps_r)|\over|\sigma_r(\eps_r)|}\right\vert
\le C_1, {\hskip1cm} 2\le j\le r,
$$ 
where $C_1=C_1(\delta)>0$ is a constant, and
$$
\left\Vert{p^s\over2\pi}\left(\mu_1\arg\sigma_j(\eps_1)
+\ldots+\mu_r\arg\sigma_j(\eps_r)\right)\right\Vert
\le C_2, {\hskip1cm} r_1\le j\le r,
$$ 
with $C_2=C_2(\delta)>0$. Set
$$
Y_j=p^s\left(\mu_1\log{|\sigma_1(\eps_1)|\over|\sigma_j(\eps_1)|}+\ldots+
\mu_r\log{|\sigma_1(\eps_r)|\over|\sigma_j(\eps_r)|}\right),
{\hskip1cm} 2\le j\le r,
$$
and
$$
Z_k={p^s\over2\pi}\left(\mu_1\arg\sigma_k(\eps_1)+\ldots+
\mu_r\arg\sigma_k(\eps_r)\right)\in\RR/\ZZ, {\hskip1cm} r_1\le k\le r.
$$
Taking $0\le\mu_i<M$, we have $M^r$ points 
$(\mu_i)_{1\le i\le r}$. The $(Y_j,Z_k)_{2\le j\le r,r_1\le k \le r}$ 
are in the product of intervals $I_j$, $2\le j\le r$,
of lengths $O(Mp^s)$ and of $r_2$ 
factors identical to $\RR/\ZZ$. 
This set can be covered by $C_3(Mp^s)^{r-1}$ 
sets of diameter at most 
$\max\{C_1,C_2 \}$, where $C_3$ is a constant
depending on $\delta$. 
By Dirichlet's {\it Schubfachprinzip}, 
choosing $M$ such that
$$
C_3(Mp^s)^{r-1}< M^r\,,
$$ 
which can be done with 
$$
M\asymp p^{(r-1)s}\,,
$$
we get that there is $(\mu_1,\ldots,\mu_r)\in\ZZ^r
\setminus \{0\}$ such that
$$
\max_{1\le i\le r}|\mu_i|\ll M\,,
$$
 $$
|Y_j|\le C_1, {\hskip1cm} 2\le j\le r,
$$ 
and
 $$
\Vert Z_k\Vert\le  C_2, {\hskip1cm} r_1\le k\le r.
$$
Set then
$$
\eta=(\eps_1^{\mu_1}\ldots\eps_r^{\mu_r})^{p^s}
$$
in such a way that $0 < \eta < 1$ (if needed,
just consider $1/\eta$).
This choice implies that
$$
|\tau_i(\eta)-1|_p=|\log_p\tau_i(\eta)|_p\le p^{-s},
\quad 0 \le i \le d,
$$
and
$$
{\left\vert\log \eta \right\vert}\ll p^s M \ll p^{rs},
$$  
and the lemma is proved. \cqfd

\proclaim Lemma 4.
Let $P(X)\in\CC_p[X]$ be a polynomial of degree $d$, and write
$$
P(X)=a_0+ \ldots +a_dX^d\,.
$$ 
Let $\xi_i$ ($1\le i\le d$) be the roots of $P(X)$ in $\CC_p$. 
Let $c$ be a real number satisfying $0 \le c \le 1$.
If
$$
|\xi_i|_p\le c\,,{\hskip1cm} 1\le i\le d,
$$ 
we get
$$
|a_k|_p\le c|a_d|_p\,,{\hskip1cm} 0\le k<d.\eqno (3.4)
$$
Conversely, if (3.4) holds, then we have
$$
|\xi_i|_p\le c^{1/d}, {\hskip1cm} 1\le i\le d.
$$

\pro
Since
$$
P(X)=a_d\prod_{1\le i\le d}(X-\xi_i)\,,
$$ 
if $|\xi_i|_p\le c\le1$ for $i=1, \ldots ,d$, then we have 
$$
|a_k|_p\le c|a_d|_p\,, \quad
\hbox{for $k=0, \ldots , d-1$}.
$$

Conversely, if 
$$
|a_k|_p\le c|a_d|_p, {\hskip1cm} 0\le k<d,
$$ 
and if $\xi\in\CC_p$ is such that 
$$
a_d\xi^d+ \ldots +a_0=0\,,
$$ 
then, there exists $k$ with
$0 \le k < d$ and
$$
|a_k\xi^k|_p\ge |a_d\xi^d|_p,
$$ 
thus,
$$
|\xi|_p^d\le|\xi|_p^{d-k}\le c\,.
$$
This completes the proof of the lemma. \cqfd

We conclude this section with two lemmas used in the
proof of Theorem 2. The first of them
was proved by Peck \cite{Peck}.

\proclaim Lemma 5.
There exists a sequence $(\eta_m)_{m \ge 1}$ of positive
units in $\cal O$ such that
$$
\eta_m\asymp e^{-dm} 
$$ 
and
$$
|\sigma_j(\eta_m)|\asymp e^m\,,{\hskip2cm}1\le j\le d\,.
$$

\pro
Let us search the unit $\eta_m$ under the form 
$$
\eta_m=\eps_1^{\mu_1}\ldots\eps_r^{\mu_r},
$$ 
with $\mu_i\in\ZZ$. 
We construct {\it real} numbers $\nu_1,\ldots,\nu_r$ such that
$$
\nu_1\log\eps_1+\ldots+\nu_r\log\eps_r=-dm\eqno (3.5)
$$ 
and
$$
\nu_1\log|\sigma_j(\eps_1)|+\ldots+\nu_r\log|\sigma_j(\eps_r)|
=m\,,{\hskip2cm} 1\le j\le d\,.\eqno (3.6)
$$ 
Taking into account that, by complex conjugation, 
the equations (3.6) corresponding to an index $j$ with $r_1\le j<r_1+r_2$ 
and to the index $j+r_2$ are identical, and that the sum of (3.5) 
and equations (3.6) is zero, we simply have to deal with a
Cramer system, since the matrix
$(\sigma_j(\eps_i))_{1\le j\le r, 1\le i\le r}$ is regular.
We solve this system and then replace every $\nu_i$ by a rational integer $\mu_i$ 
such that $|\mu_i-\nu_i|\le1/2$. \cqfd

\proclaim Lemma 6.
Let $\lambda'$ be a positive real number.
Let $(\eta_m)_{m \ge 1}$ be a sequence of positive units 
as in Lemma 5.
There exists a finite set $\Gamma=\Gamma(\lambda')$ 
of non-zero elements of $\K$ 
such that for all integer polynomial $P(X)$ 
of degree at most $d$ that satisfies 
$$
|P(\alpha)| \le \lambda' H(P)^{-d},  \eqno (3.7)
$$
there
exist a positive integer $m$ and $\gamma$ in $\Gamma$ for which
$$
P(\alpha)=\gamma\eta_m\,.
$$ 

\pro 
Below, all the constants implicit in $\ll$ depend on $\K$ and on $\lambda'$.
Let $m$ be a positive integer such that
$$
H(P)\asymp e^m\,,
$$ 
and set
$$
\gamma=P(\alpha)\eta_m^{-1}.
$$ 
Since $D\alpha^k$ is an algebraic integer for $k=0,\ldots,d$, 
the algebraic number $D\gamma$ 
is an algebraic integer, and, by (3.7), 
$$
|\gamma|\ll1\,.
$$
Furthermore, for $j = 1, \ldots ,  d$, we have
$$
|\sigma_j(\gamma)|=|P(\sigma_j(\alpha))| \cdot
|\sigma_j(\eta_m^{-1})|\ll H(P)e^{-m}\ll1\,.
$$ 
The algebraic integers $D\gamma\in\cal O$
and all their complex conjugates being bounded, 
they form a finite set. \cqfd

\vskip 6mm

\centerline{\bf 4. Proof of Theorem 1}

\vskip 4mm

Let $\delta$ be in $(0, 1)$ to be selected later.
Apply Lemma 3 with this $\delta$ to get a unit
$\eta$ and apply Lemma 1 with this unit and with
$\gamma = \alpha^d$.
Since $D^2 \alpha^d \eta \in \ZZ + \ldots + \alpha^d \ZZ$,
we get
$$
D^2\eta\alpha^d=a_0+a_1\alpha+ \ldots +a_d\alpha^d
=P(\alpha),
$$ 
where, by Corollary 1, $P(X)$ is an integer polynomial of degree $d$ and
$$
|P(\alpha)| \asymp H(P)^{-d} \asymp H^{-d}.
$$
By Lemma 2, each coefficient of $P(X)$
has its $p$-adic
absolute value $\ll (\log 3 H(P))^{-1/r}$,
except the leading coefficient, whose
$p$-adic absolute value equals $|D|_p^2$.

We then infer from Lemma 4 that all the roots
of $P(X)$ in $\CC_p$ have their
$p$-adic absolute value $\ll (\log 3 H(P))^{-1/(dr)}$.
This proves (2.2).

It now remains for us to 
guarantee that $P(X)$ has a root very close to
$\alpha$. To this end, we proceed to check that
$$
|P'(\alpha)|\gg H(P)\,.
$$ 
Since
$$
P'(\alpha)=a_1+ \ldots +da_d\alpha^{d-1}\,,
$$ 
we get
$$
P'(\alpha)=D^2\left({\rm Tr}(\eta\alpha^d\beta_1)+
2\alpha{\rm Tr}(\eta\alpha^d\beta_2)
+ \ldots +d\alpha^{d-1}{\rm Tr}(\eta\alpha^d\beta_d)\right),
$$ 
hence,
$$
P'(\alpha)=D^2\sum_{i=0}^d\sum_{k=1}^dk\alpha^{k-1}\sigma_i(\eta\alpha^d\beta_k)\,.
$$
Let us write 
$$
P'(\alpha)=D^2\sum_{i=0}^dA_i\sigma_i(\eta)
$$ 
with
$$
A_i=\sigma_i(\alpha^d)\sum_{k=1}^dk\alpha^{k-1}\sigma_i(\beta_k)\,,
\quad i = 0, \ldots , d.
$$
Observe first that
$$
\sum_{i=1}^dA_i\not=0\,.
$$
Indeed, if this is not the case, then, working with the 
unit $\eta=1$, that is, 
with $P(X)=D^2X^d$ and $P'(\alpha)=dD^2\alpha^{d-1}$, 
we get 
$$
d\alpha^{d-1}=A_0=\alpha^d\sum_{k=1}^dk\alpha^{k-1}\beta_k,
$$ 
hence,
$$
d=\sum_{k=1}^dk\alpha^k\beta_k\,.
$$ 
Taking the trace, and recalling that ${\rm Tr}(\alpha^k\beta_k)=1$, we get 
$$
d(d+1)=\sum_{k=1}^dk\,,
$$ 
a contradiction.

Write
$$
P'(\alpha)=D^2\sum_{i=1}^dA_i\sigma_i(\eta)+O(H^{-d})=
D^2\sigma_1(\eta)\sum_{i=1}^dA_i+B
$$ 
with
$$
|B|\le D^2\sum_{2\le i\le d}|A_i| \cdot |\sigma_1(\eta)| \cdot
\biggl|{\sigma_i(\eta) \over \sigma_1(\eta)} - 1 \biggr|
+O(H^{-d})\,.
$$ 
Selecting now $\delta$ such that
$$
\delta\sum_{2\le i\le d}|A_i| \le
{1\over3}\left\vert\sum_{i=1}^dA_i\right\vert,
$$ 
we infer from Lemma 3 that
$$
|P'(\alpha)|\ge{1\over 2}D^2\left\vert\sigma_1(\eta)\sum_{i=1}^dA_i\right\vert,
$$ 
when $H$ is sufficiently large.
This gives
$$
|P'(\alpha)|\gg |\sigma_1(\eta)| \gg H\,.
$$
Consequently, $P(X)$ has a root $\xi$ such that 
$$
|\alpha-\xi|\ll H(P)^{-d-1}\ll H(\xi)^{-d-1}.
$$
Classical arguments (see at the end of the proof of Theorem 2.11
in \cite{BuLiv}) show that $\xi$ must be real 
and of degree $d$ if $H$ is sufficiently
large.
This proves (2.1). Inequality (2.3) follows from
(2.1) and (2.2) together with the fact that $\xi$
is of degree $d$. This completes the proof of the theorem.

\vskip 6mm

\centerline{\bf 5. Proof of Theorem 2}

\vskip 4mm

The constants implicit in $\ll$ and $\gg$
below depend on $\K$, $p$ and $\lambda$.
There exists a positive real number $\lambda'$,
depending on $\lambda$ and on $d$, such that
the minimal defining polynomial $P(X)$ of any real number
$\xi$ of sufficiently large height and 
for which (2.4) holds
is of degree $d$ and satisfies
$$
|P(\alpha)| \le \lambda' H(P)^{-d}.
$$
Let $(\eta_m)_{m \ge 1}$ be as in Lemma 5. By Lemma 6,
it is sufficient to prove Theorem 2 for the integer polynomials 
$P(X)$ as above such that
$$
P(\alpha)=\gamma\eta_m=a_0+a_1\alpha+\ldots+a_d\alpha^d.
$$ 
Let $\xi_i$ be the roots of $P(X)$ in $\CC_p$ and set
$$
u := \max_{1 \le i \le d} \, |\xi_i|_p\,.
$$ 
Assume that $u\le1\,.$ It follows from Lemma 4 that 
$$
|a_k|_p\le u|a_d|_p\,,{\hskip2cm}0\le k< d\,.
$$ 
Dividing $P(X)$ by $p^s=|a_d|_p^{-1}$ if necessary, 
we can assume that $|a_d|_p=1$, and we obtain that
$$
|a_k|_p\le u\,,{\hskip2cm}0\le k< d \,.
$$ 
For $j = 1, \ldots , d$, we then have 
$$
\gamma\eta_m\alpha^{-d}-\tau_j(\gamma\eta_m\alpha^{-d})=
\sum_{k=0}^{d-1}a_k(\alpha^{k-d}-\tau_j(\alpha^{k-d})),
$$ 
hence,
$$
|\gamma\eta_m\alpha^{-d}-\tau_j(\gamma\eta_m\alpha^{-d})|_p\ll u\,.
$$ 
Since $|\eta_m|_p=1$, we get that
$$
\left\vert{\tau_j(\eta_m)\over\eta_m}{\tau_j(\gamma)\alpha^d
\over\gamma\tau_j(\alpha^d)}-1\right\vert_p\ll u\,.
$$ 
Upon writing  
$$
\eta_m=\eps_1^{\mu_{1,m}}\ldots\eps_r^{\mu_{r,m}}\,,
$$ 
we have thus 
$$
u\gg\left\vert\left({\tau_j(\eps_1)\over\eps_1}\right)^{-\mu_{1,m}}
\ldots\left({\tau_j(\eps_r)\over\eps_r}\right)^{-\mu_{r,m}}
{\tau_j(\gamma)\alpha^d\over\gamma\tau_j(\alpha^d)}-1\right\vert_p
$$ 
If 
$$
{\tau_j(\eta_m)\over\eta_m}=
{\gamma\tau_j(\alpha^d)\over\tau_j(\gamma)\alpha^d}
$$ 
holds for $j=1,\ldots,d$, the number 
$$
\gamma\eta_m\alpha^{-d}
$$ 
is equal to all its conjuguates, hence is rational, 
and we have 
$$
P(\alpha)=b\alpha^d
$$ 
with $b\in\QQ$, hence $P(X)=bX^d$, a contradiction. For every $m$, 
there thus exists an index $j$ such that $1\le j\le d$ and 
$$
\left({\tau_j(\eps_1)\over\eps_1}\right)^{-\mu_{1,m}}\ldots
\left({\tau_j(\eps_r)\over\eps_r}\right)^{-\mu_{r,m}}
{\tau_j(\gamma)\alpha^d\over\gamma\tau_j(\alpha^d)}\not=1\,.
$$ 
Consequently, by the theory of linear forms in non-Archimedean
logarithms (see e.g., Kunrui Yu \cite{Yu}), there exists a positive constant 
$\kappa$ such that 
$$
u\gg(\max_{1\le i\le r}|\mu_{i,m}|)^{-\kappa}.
$$ 
Since $\eta_m\asymp H(P)^{-d}$ and 
$$
|\log\eta_m|\asymp\max_{1\le i\le r}|\mu_{i,m}|,
$$ 
the matrix $(\log|\sigma_j(\eps_i)|)_{1\le i\le r, 1\le j\le r}$ being regular, 
we conclude that 
$$
u\gg(\log 3 H(\xi))^{-\kappa}\,.
$$
This completes the proof of the theorem. \cqfd


\vskip 7mm

\centerline{\bf References}

\vskip 7mm

\beginthebibliography{999}

\bibitem{Bos94}
M. D. Boshernitzan, 
{\it Elementary proof of Furstenberg's Diophantine result}, 
Proc. Amer. Math. Soc.  122  (1994),  67--70.

\sm

\bibitem{BuLiv}
Y. Bugeaud,
Approximation by algebraic numbers,
Cambridge Tracts in Mathematics 160,
Cambridge, 2004.

\sm

\bibitem{BDM}
Y. Bugeaud, M. Drmota, and B. de Mathan,
{\it On a mixed Littlewood conjecture in Diophantine approximation},
Acta Arith. 128 (2007), 107--124.

\sm

\bibitem{EiKl07}
M. Einsiedler and D. Kleinbock,
{\it Measure rigidity and $p$-adic Littlewood-type problems},
Compositio Math. 143 (2007), 689--702.

\sm

\bibitem{Fur67}
H. Furstenberg,
{\it Disjointness in ergodic theory, minimal sets, and a problem
in Diophantine approximation},
Math. Systems Theory 1 (1967), 1--49.

\sm

\bibitem{BdM05}
B. de Mathan,
{\it On a mixed Littlewood conjecture for quadratic numbers},
J. Th\'eor. Nombres Bordeaux 17 (2005), 207--215.

\sm

\bibitem{BdMTe}
B. de Mathan et O. Teuli\'e,
{\it Probl\`emes diophantiens simultan\'es},
Monatsh. Math. 143 (2004), 229--245.

\sm

\bibitem{Peck}
L. G. Peck,
{\it Simultaneous rational approximations to algebraic numbers},
Bull. Amer. Math. Soc. 67 (1961), 197--201. 

\sm

\bibitem{Yu}
Kunrui Yu,
{\it $p$-adic logarithmic forms and group varieties II}, 
Acta Arith. 89 (1999), 337--378. 

\sm

\bibitem{Wir}
{E. Wirsing},
{\it Approximation mit algebraischen Zahlen beschr\"ankten
Grades}, J. reine angew. Math. {206} (1961), 67--77.

\endthebibliography

\goodbreak

\vskip 8mm

Yann Bugeaud \hfill Bernard de Mathan

Universit\'e Louis Pasteur \hfill Universit\'e Bordeaux I

Math\'ematiques \hfill Institut de Math\'ematiques

7, rue Ren\'e Descartes \hfill 351, cours de la Lib\'eration

67084 STRASBOURG Cedex (France) \hfill33405 TALENCE Cedex (France)

\medskip

{\tt bugeaud@math.u-strasbg.fr} \hfill {\tt demathan@math.u-bordeaux1.fr}

\bye